\documentstyle[amssymb]{article}

\begin{document}
\font\frak=eufm10 scaled\magstep1
\font\fak=eufm10 scaled\magstep2
\font\fk=eufm10 scaled\magstep3
\font\scriptfrak=eufm10
\font\tenfrak=eufm10

\newtheorem{theorem}{Theorem}
\newtheorem{corollary}{Corollary}
\newtheorem{proposition}{Proposition}
\newtheorem{definition}{Definition}
\newtheorem{lemma}{Lemma}
\font\frak=eufm10 scaled\magstep1
\newenvironment{pf}{{\noindent{\it Proof. }}}{\ $\Box$\medskip}


\mathchardef\za="710B  
\mathchardef\zb="710C  
\mathchardef\zg="710D  
\mathchardef\zd="710E  
\mathchardef\zve="710F 
\mathchardef\zz="7110  
\mathchardef\zh="7111  
\mathchardef\zvy="7112 
\mathchardef\zi="7113  
\mathchardef\zk="7114  
\mathchardef\zl="7115  
\mathchardef\zm="7116  
\mathchardef\zn="7117  
\mathchardef\zx="7118  
\mathchardef\zp="7119  
\mathchardef\zr="711A  
\mathchardef\zs="711B  
\mathchardef\zt="711C  
\mathchardef\zu="711D  
\mathchardef\zvf="711E 
\mathchardef\zq="711F  
\mathchardef\zc="7120  
\mathchardef\zw="7121  
\mathchardef\ze="7122  
\mathchardef\zy="7123  
\mathchardef\zf="7124  
\mathchardef\zvr="7125 
\mathchardef\zvs="7126 
\mathchardef\zf="7127  
\mathchardef\zG="7000  
\mathchardef\zD="7001  
\mathchardef\zY="7002  
\mathchardef\zL="7003  
\mathchardef\zX="7004  
\mathchardef\zP="7005  
\mathchardef\zS="7006  
\mathchardef\zU="7007  
\mathchardef\zF="7008  
\mathchardef\zW="700A  

\newcommand{\be}{\begin{equation}}
\newcommand{\ee}{\end{equation}}
\newcommand{\ra}{\rightarrow}
\newcommand{\lra}{\longrightarrow}
\newcommand{\bea}{\begin{eqnarray}}
\newcommand{\eea}{\end{eqnarray}}
\newcommand{\beas}{\begin{eqnarray*}}
\newcommand{\eeas}{\end{eqnarray*}}
\newcommand{\Z}{{\Bbb Z}}
\newcommand{\R}{{\Bbb R}}
\newcommand{\C}{{\Bbb C}}
\newcommand{\1}{{\bold 1}}
\newcommand{\SL}{SL(2,\C)}
\newcommand{\Sl}{sl(2,\C)}
\newcommand{\SU}{SU(2)}
\newcommand{\su}{su(2)}
\newcommand{\G}{{\goth g}}
\newcommand{\D}{{\rm d}}
\newcommand{\de}{\,{\stackrel{\rm def}{=}}\,}
\newcommand{\we}{\wedge}
\newcommand{\nn}{\nonumber}
\newcommand{\ot}{\otimes}
\newcommand{\s}{{\textstyle *}}
\newcommand{\ts}{T^\s}
\newcommand{\da}{\dagger}
\newcommand{\pa}{\partial}
\newcommand{\ti}{\times}
\newcommand{\A}{{\cal A}}
\newcommand{\E}{{\cal E}}
\newcommand{\Li}{{\cal L}}
\newcommand{\ka}{{\Bbb K}}
\newcommand{\find}{\mid}
\newcommand{\ad}{{\rm ad}}
\newcommand{\rS}{]^{SN}}
\newcommand{\rb}{\}_P}
\newcommand{\p}{{\sf P}}
\newcommand{\h}{{\sf H}}
\newcommand{\X}{{\cal X}}
\newcommand{\I}{\,{\rm i}\,}
\newcommand{\rB}{]_P}
\newcommand{\qd}{{\rm Qder}_\A(\E)}
\newcommand{\di}{{\rm Diff}^1_\A(\E)}
\title{Quasi-derivations and QD-algebroids}
\author{Janusz Grabowski\thanks{Supported by KBN, grant No 2 PO3A 041 18}\\
Mathematical Institute, Polish Academy of Sciences\\
ul. \'Sniadeckich 8, P.O. Box 21, 00-956 Warszawa 10, Poland\\
{\it e-mail:} jagrab@impan.gov.pl}
\date{}
\maketitle
\begin{abstract}\noindent
Axioms of Lie algebroid are discussed. In particular, it is shown
that a Lie  QD-algebroid (i.e. a  Lie   algebra bracket on the
$C^\infty(M)$-module $\E$ of sections of a  vector  bundle $E$
over  a manifold $M$ which satisfies $[X,fY]=f[X,Y]+A(X,f)Y$ for
all $X,Y\in\E$,   $f\in C^\infty(M)$, and for certain $A(X,f)\in
C^\infty(M)$) is a Lie algebroid if $rank(E)>1$, and is a local
Lie  algebra in  the sense of Kirillov if $E$ is a line bundle.
Under a weak condition also the skew-symmetry of the bracket is
relaxed.

\end{abstract}
\section{Introduction}
A {\it Lie algebroid}  (or  its  pure  algebraic   counterpart --
a {\it   Lie pseudoalgebra}) is an object so natural that the
people often use it not even realizing this fact.

\bigskip\noindent
{\bf Definition.} Let $R$ be a commutative and unitary ring, and
let $\A$ be a commutative  $R$-algebra. A {\it Lie pseudoalgebra}
over $R$ and $\A$ is an $\A$-module $\E$ together with a bracket
$[\cdot,\cdot] :\E\ti\E\ra\E$ on the module $\E$, and an
$\A$-module  morphism $a:\E\ra{\rm Der}(\A)$ from $\E$ to the
$\A$-module ${\rm Der}(\A)$ of derivations of $\A$, called the
{\it anchor} of $\E$, such that \noindent
\begin{description}
\item{(i)} the bracket on $\E$ is $R$-bilinear, alternating, and
satisfies the Jacobi identity:
\be\label{ji}
[[X,Y],Z]=[X,[Y,Z]]-[Y,[X,Z]].
\ee
\item{(ii)} For all  $X,Y\in  \E$  and  all $f\in\A$ we have
\be\label{0}
[X,fY]=f[X,Y]+a(X)(f)Y;
\ee
\item{(iii)}   $a([X,Y])=[a(X),a(Y)]_c$   for   all    $X,Y\in\E$,    where
$[\cdot,\cdot]_c$ is the commutator bracket on ${\rm Der}(\A)$.
\end{description}
A {\it Lie algebroid} on a vector bundle $E$ over a base manifold
$M$ is  a Lie pseudoalgebra over the $\R$-algebra of smooth
functions $\A=C^\infty(M)$ on the $C^\infty(M)$-module $\E=Sec(E)$
of smooth sections of $E$. In this case the  anchor map  is
described  by  a  vector  bundle morphism $a:E\ra TM$ which
induces the bracket homomorphism from $(\E, [\cdot,\cdot])$ into
the Lie algebra $(\X(M),[\cdot,\cdot]_{vf})$ of vector fields on
$M$.

\smallskip
Lie pseudoalgebras appeared first in the paper of Herz  \cite{He} but
one can find similar concepts under more than a dozen  of names  in
the   literature   (e.g. Lie modules,  $(R,A)$-Lie   algebras,
Lie-Cartan pairs, Lie-Rinehart algebras, differential algebras, etc.).
Lie algebroids were introduced  by Pradines \cite{Pr} as infinitesimal
parts of differentiable groupoids. In the same year a book by Nelson
was published where a general theory of Lie modules together with a
big part of the corresponding differential calculus can be found. We
also refer to a survey article by Mackenzie \cite{Ma}.

Lie algebroids on a singleton base space are Lie algebras. Another
canonical example is the tangent bundle $TM$ with the canonical
bracket $[\cdot,\cdot]_{vf}$ on the space $\X(M)=Sec(TM)$ of
vector fields.

\smallskip During  Poisson seminars at the Warsaw University we
realized long ago that the axiom (iii) telling that the anchor map
induces a Lie algebra homomorphism is superfluous for Lie
algebroids and that also the rest of assumptions can  be weakened.
The proofs are very elementary and we  did  not pay much attention
to this result, especially when we discovered that such
observations, although not explicitly in the language of Lie
algebroids, had been already done in the old paper by Herz [He]
and repeated by Kosmann-Schwarzbach and Magri ([KSM], section
6.1). However, to our surprise, all the time we find some people,
even people  working in the subject, who seem to have no knowledge
about  these  simple  facts, or rediscover them, or, in the worst
case, try to build generalized object which are impossible due to
these observations. The fact that these remarks concerning the
axioms  of  a  Lie algebroid seem not to be commonly known,
encouraged us finally  to  publish the  present  note (circulating
since a few years as a private communication) in which we  discuss
axioms of Lie algebroid and present a survey of relevant  results
in this direction. Our setting is slightly more general than  just
Lie algebroids, since  we accept also certain brackets which do
not satisfy the Jacobi identity  (what is already present in
[KSM]) or are not assumed to  be skew-symmetric.  We add also
certain results on the axioms of Poisson and Jacobi brackets which
can be of some interest in this context.

\section{Quasi-derivations in modules and vector bundles}
Let $\A$ be an associative ring (not necessarily commutative) over a
commutative and unitary ring  $R$ and let $\E$ be an $\A$-module. Our
standard example will be the $C^\infty(M)$-module of sections of a
vector bundle $E$ over $M$. We will write the operation in $\A$ and
the $\A$-action on $\E$ without special signs: $\A\ti\A\ni(f,g)\mapsto
fg\in\A$ and $\A\ti\E\ni(f,X)\mapsto fX\in\E$. For  $f\in\A$ we denote
by $f_\E$ the operator on $\E$ afforded by the module structure:
$f_\E(X)=fX$. Let $\A_\E$ denote the set of all such operators. We
call a linear operator $D:\E\ra\E$  a {\it quasi-$\A$-derivation} if
$[D,\A_\E]_c\subset\A_\E$, where the bracket operation
$[\cdot,\cdot]_c$ is the commutator. Let $\qd$ denote the $R$-module
of all quasi-$\A$-derivations. We have the following (cf. \cite{Sk},
Proposition 1.1).

\begin{theorem} If $\E$ is a  faithful  $\A$-module,  then  there  is  a
linear mapping $\qd\ni D\mapsto\widehat{D}\in {\rm Der}(\A)$ such
that
$$ [D,f_\E]_c=(\widehat{D}(f))_\E.
$$
\end{theorem}
\begin{pf} By definition, for any $D\in\qd$ there is $\widehat{D}(f)\in\A$
such that $[D,f_\E]_c=\widehat{D}(f)_\E$ and $\widehat{D}(f)$ is
uniquely determined for any $f\in\A$,  since  the  module  is
faithful. We have  $D(fgX)=fgD(X)+\widehat{D}(fg)X$. On the other
hand,
$$
D(fgX)=fD(gX)+\widehat{D}(f)gX=fgD(X)+(f\widehat{D}(g)+\widehat{D}(f)g)X,
$$
so that $\widehat{D}(fg)=f\widehat{D}(g)+\widehat{D}(f)g$, since the module
is faithful, i.e. $\widehat{D}$ is an $\A$-derivative. The mapping
$D\mapsto  \widehat{D}$  is clearly linear.
\end{pf}

Let us assume for the future that we deal with faithful modules. The
above map $\qd\ni D\mapsto\widehat{D}\in {\rm Der}(\A)$ we will call
the {\it universal anchor map} for the $\A$-module $\E$. The universal
anchor map induces the short exact sequence
$$
0\ra{\rm End}_\A(\E)\ra\qd\ra{\rm Der}(\A).
$$
The following theorem shows that this is  a  short  exact  sequence  of  Lie
algebras.

\begin{theorem}\label{t1} The $R$-module $\qd$ is closed with respect to the
commutator. Moreover, the universal anchor map is a Lie algebra homomorphism,
i.e.
$$
\widehat{[D_1,D_2]_c}=[\widehat{D_1},\widehat{D_2}]_c
$$
for $D_1,D_2\in\qd$.
\end{theorem}
\begin{pf} By definition,
\beas
[D_1,D_2]_c(fX)&=&D_1(D_2(fX))-D_2(D_1(fX))\\
&=&D_1(fD_2(X)+\widehat{D_2}(f)X)-D_2(fD_1(X)+\widehat{D_1}(f)X)\\
&=&f(D_1(D_2(X))-D_2(D_1(X))+(\widehat{D_1}(\widehat{D_2}(f))-
\widehat{D_2}(\widehat{D_1}(f)))X\\
&=&f[D_1,D_2]_c(X)+[\widehat{D_1},\widehat{D_2}]_c(f)X.
\eeas
Hence, $[D_1,D_2]_c$ is a quasi-derivation and
$$
\widehat{[D_1,D_2]_c}(f)X=[\widehat{D_1},\widehat{D_2}]_c(f)X
$$
for all $X\in\E$ and all $f\in\A$ and the theorem follows, since the
module is faithful.
\end{pf}

In  case  $\A$  is  a  commutative ring,   $\qd$  is  additionally   an
$\A$-module with respect  to the  action $fD=f_\E\circ D$.

\begin{proposition} If $\A$  is  a  commutative  ring  then  the  commutator
bracket on the $\A$-module $\qd$ satisfies
$$
[D_1,fD_2]_c=f[D_1,D_2]_c+\widehat{D_1}(f)D_2.
$$
\end{proposition}
\begin{pf} It follows from the Leibniz rule for the commutator bracket:
$$
[D_1,fD_2]_c=[D_1,f_\E\circ D_2]_c=
[D_1,f_\E]_c\circ D_2+f_\E\circ[D_1,D_2]_c. $$
\end{pf}

\begin{corollary} If $\A$ is commutative then $\qd$ is  canonically  a  Lie
pseudoalgebra over $R$ and $\A$.
\end{corollary}
It is pointed out in \cite{KM} that concept of quasi-derivation can be
traced back to N.~Jacobson \cite{Ja1,Ja2} as a special case of his
{\it pseudo-linear endomorphism}. It has appeared also in \cite{Ne}
under the name of a {\it module derivation} and used to define linear
connections in the algebraic setting. In the geometric setting for Lie
algebroids it has been studied in \cite{Mk}, Ch. III, under the name
{\it covariant differential operator}. For more detailed history and
recent development we refer to \cite{KM} where quasi-derivations are
called {\it derivative endomorphisms}.

\section{QD-algebroids}
Suppose now that $\E$ is faithful and that on $\E$ we have
additionally a bracket operation  $[\cdot,\cdot]$  which  is
$R$-linear  in  the  second argument and satisfies the Jacobi
identity (\ref{ji}), i.e. $\ad_{[X,Y]}=[\ad_X,\ad_Y]_c$, where
$\ad_X(Y)=[X,Y]$ is  the  adjoint operator. Let us assume that the
adjoint operators $\ad_X$ are in $\qd$ for any $X\in\E$. We can
define the {\it anchor  map} $\E\ni X\mapsto\widehat{X}\in{\rm
Der}(\A)$ by $\widehat{X}=\widehat{\ad_X}$.
\begin{proposition}\label{p1} The anchor map is a bracket homomorphism, i.e.
$$
\widehat{[X,Y]}=[\widehat{X},\widehat{Y}]_c.
$$
for all $X,Y\in\E$.
\end{proposition}
\begin{pf} This is a direct consequence of Theorem \ref{t1} and the  Jacobi
identity.
\end{pf}

Consider a particular example, where $\A$ is the commutative
associative algebra  of  smooth  functions  on  a  manifold  $M$ and
$\E$  is  the $\A$-module  of  smooth  sections  of  a vector bundle
$E$  over  $M$ with $k$-dimensional fibers.  Note that this  module is
faithful. In this case quasi-derivations are particular first-order
differential operators on $E$.  An $\R$-bilinear bracket
$[\cdot,\cdot]$ on  $\E$  we   will call   a {\it quasi-derivation
algebroid} (shortly {\it QD-algebroid})  if   the bracket is a
quasi-$\A$-derivation of $\E$ for both variables separately. Hence,
$[X,fY]=f[X,Y]+\widehat{X}(f)Y$    and
$[fX,Y]=f[X,Y]+\widetilde{Y}(f)X$ for   all  $X,Y\in\E$  and  all
$f\in\A$, where $X\mapsto\widehat{X}\in{\rm Der}(\A)$ and
$Y\mapsto\widetilde{Y}\in{\rm Der}(\A)$ are, respectively, the left
and the right anchor maps.  In  the  case  when  the  anchor maps are
tensorial ($\A$-linear),  we  speak  about an {\it algebroid} (cf.
\cite{GU}).
\par
In the case when the bracket is a Lie algebra bracket, we will
speak about  a {\it Lie QD-algebroid} (resp., {\it Lie
algebroid}). We can also consider {\it Loday QD-algebroid} (resp.,
{\it  Loday  algebroid})  requiring  additionally only the Jacobi
identity  (\ref{ji}) without skew-symmetry assumption. The name
is   motivated by  {\it Loday algebras}, i.e. the concept of
`non-skew-symmetric Lie algebras' introduced by J.~L.~Loday
\cite{Lo1,Lo2} (and called by Loday {\it Leibniz algebras}).

\begin{theorem} Every QD-algebroid (resp., Loday QD-algebroid, or  Lie
QD-algebroid) of rank $>1$ is an algebroid (resp., Loday
algebroid, or Lie algebroid).
\end{theorem}
\begin{pf} According to Theorem \ref{t1}, we have
$[X,fY]=f[X,Y]+\widehat{X}(f)Y$  and
$[fX,Y]=f[X,Y]+\widetilde{Y}(f)X$ for  all $X,Y\in\E$ and all
$f\in\A$, where $\widehat{X}$ and $\widetilde{Y}$ are derivations
of $\A$. Moreover, due to Proposition \ref{p1}, in the case of a
Loday QD-algebroid we have
$\widehat{[X,Y]}=[\widehat{X},\widehat{Y}]_c$. For a Lie
QD-algebroid, additionally, $\widehat{X}=-\widetilde{X}$. It
remains to   show   that   for   QD-algebroids   of   rank $>1$
the    maps $X\mapsto\widehat{X}$ and $Y\mapsto \widetilde{Y}$ are
$\A$-linear.   To see  this,  let  us  write $[gX,fY]$ in two
different ways. First,
$$
[gX,fY]=f[gX,Y]+\widehat{gX}(f)Y=fg[X,Y]+f\widetilde{Y}(g)X+\widehat{gX}(f)Y.
$$
On the other hand,
$$
[gX,fY]=g[X,fY]+\widetilde{fY}(g)X=gf[X,Y]+g\widehat{X}(f)Y+
\widetilde{fY}(g)X,
$$
so that
$$
f\widetilde{Y}(g)X+\widehat{gX}(f)Y=g\widehat{X}(f)Y+\widetilde{fY}(g)X
$$
or \be\label{1}
(\widehat{gX}-g\widehat{X})(f)Y=(\widetilde{fY}-f\widetilde{Y})(g)X.
\ee But, since the rank of the bundle is $>1$, for  a  given
$X\in\E$  we  can choose $Y$ which is locally linearly independent
(and vice versa), so  that (\ref{1})      implies
$\widehat{gX}-g\widehat{X}=0$        and
$\widetilde{fY}-f\widetilde{Y}=0$ for any $X,Y\in\E$ and $f\in\A$.
\medskip
\end{pf}

As it has been proved in \cite{GM1}, the difference between Loday algebroids
and Lie algebroids is not very big.

\begin{theorem} For any Loday algebroid, the  left  and  the  right  anchor
differ  by  sign,  $\widetilde{X}=-\widehat{X}$,   and   the   bracket   is
skew-symmetric at points of $M$ where they do not vanish.
\end{theorem}

In other words, a Loday algebroid is, in  fact  a  Lie algebroid
around points where the anchors (or, equivalently, one of the
anchors)  do not vanish, and it has Loday  algebras  as fibers
over points  where  the anchors do vanish.

\medskip
Consider now QD-algebroids of rank 1. The bracket defined on
sections of a line-bundle $E$ over $M$ is a differential operator
of first order, so local. Assuming the Jacobi identity for  the
bracket,  i.e.  dealing  with Loday QD-algebroid, one can prove,
in fact, skew-symmetry of the bracket (\cite{GM2}, Corollary 2).

\begin{theorem}\label{t5}  Every  Loday  QD-algebroid  of  rank   1   is
a   Lie QD-algebroid.
\end{theorem}

A Lie QD-algebroid of rank 1 is just a local Lie algebra in the
sense  of  Kirillov \cite{Ki}. In the trivial case $E=M\ti\R$ the
bracket is defined  on $C^\infty(M)=Sec(E)$ and it is called a
Jacobi bracket \cite{Li}. It is of  the  form \be\label{jb}
[f,g]=\zL(\D  f,\D  g)+f\zG(g)-g\zG(f) \ee for  certain bivector
field $\zL$ and certain vector field $\zG$. The Jacobi identity is
equivalent to compatibility conditions, usually  written  in terms
of  the Schouten-Nijenhuis  bracket:
$$
[\zG,\zL]^{SN}=0, \quad [\zL,\zL]^{SN}+2\zL\we\zG=0.
$$
One has $\widehat{f}=i_{\D f}\zL+f\zG$. In the case when $\zG=0$
we deal with  a {\it Poisson bracket} and in the case $\zL=0$ --
with a Lie algebroid bracket. This means that Lie algebroids on
the trivial 1-dimensional bundle (i.e. Lie algebroid brackets on
$C^\infty(M)$) are of the form $[f,g]=f\zG(g)-g\zG(f)$ for certain
vector field $\zG$ on $M$. Here the anchor map is
$\widehat{f}=f\zG$. Other natural examples of rank-1 Lie
algebroids are associated with foliations of $M$ by
one-dimensional leaves.

Of course, Theorem \ref{t5} implies in particular that we can skip
skew-symmetry in the definition of Jacobi and Poisson brackets.
The direct proof for Poisson brackets is particularly simple.
Indeed, the Jacobi identity (\ref{ji}) implies immediately that
$[[f,g]+[g,f],h]=0$ for all $f,g,h\in C^\infty(M)$. Using the
Leibniz rule with respect to both arguments, we get
$$0=[[f^2,g]+[g,f^2],h]=2([f,g]+[g,f])[f,h].
$$
For $f=g=h$ we get $4[f,f]^3=0$, thus $[f,f]=0$ and skew-symmetry
follows.

\smallskip
Let us summarize the above results as follows.

\begin{theorem} Let $E$ be a vector bundle over $M$ and let $[\cdot,\cdot]$
be a bilinear bracket operation  on  the  $C^\infty(M)$-module  $\E=Sec(E)$
that satisfies the Jacobi identity (\ref{ji}) and that
is a quasi-derivation with respect to both arguments.
\begin{description}
\item{(a)} If $rank(E)>1$, then there is a vector bundle  morphism  $a:E\ra
TM$ over the identity map on  $M$  such  that  $a([X,Y])=[a(X),a(Y)]_{vf}$
and
$$
[fX,gY]=fg[X,Y]+fa(X)(g)Y-ga(Y)(f)X
$$
for all $X,Y\in\E$, $f,g\in C^\infty(M)$. Moreover, $[X,Y](p)=-[Y,X](p)$ if
$a_p\ne 0$.
\item{(b)} If $rank(E)=1$, then the bracket is skew-symmetric and defines a
local Lie algebra structure which,  locally,  is  equivalent  to  a  Jacobi
bracket (\ref{jb}).
\end{description}
\end{theorem}
\begin{corollary}
Lie QD-algebroids on $E$ are exactly Lie algebroids if
$rank(E)>1$, and local  Lie  algebras  in  the  sense  of Kirillov
if  $rank(E)=1$.
\end{corollary}
\begin{corollary}
A Lie algebroid on a vector bundle  $E$ of rank $>1$ is just a Lie
bracket on sections of $E$ which is a quasi-derivation with
respect to one (hence both) argument.
\end{corollary}
\medskip\noindent
{\bf Acknowledgments.} The author wish to thank the referees for
useful comments.

\end{document}